\documentclass[a4paper,11pt]{article}

\usepackage{amssymb}

\newtheorem{theoreme}{Theorem}[section]

\newenvironment{demo}{\noindent {\sl Proof}. \ }{\qed}

\newenvironment{ex}{\noindent {\bf Example}. \ }{\\}
\newtheorem{stheoreme}{Theorem}[subsection]
\newtheorem{sdefin}[stheoreme]{Definition}
\newtheorem{sprop}[stheoreme]{Proposition}
\newtheorem{slemme}[stheoreme]{Lemma}
\newtheorem{scorol}[stheoreme]{Corollary}

\font\tenCal=cmsy10
\newfam\Calfam
\textfont\Calfam=\tenCal

\def\Id{\mathop{\hbox{\textrm Id}}\nolimits}
\def\Aut{\mathop{\hbox{\textrm Aut}}\nolimits}

\def\Frac{\mathop{\hbox{\textrm Frac}}\nolimits}

\def\qed{\hfill{$\sqcap\!\!\!\!\sqcup$}}

\def\a{\alpha}
\def\b{\beta}

\title {Isomorphisms between quantum generalized Weyl
algebras}

\author{Lionel Richard\footnote{Institut Girard Desargues,
Universit\'e Lyon 1.
F-69622 Villeurbanne Cedex - France.
richard@igd.univ-lyon1.fr} \
and 
Andrea Solotar\footnote{Dto. de Matem\'atica, 
Facultad de Cs. Exactas y Naturales.
Universidad de Buenos Aires.
Ciudad Universitaria Pab I. (1428),
Buenos Aires - Argentina.
asolotar@dm.uba.ar}}

\begin{document}\maketitle
\begin{footnotetext}{
{Research partially supported 
by \textsc{UBACyT TW62}, \textsc{PICS-CNRS 1514}, \textsc{Fundaci\'on Antorchas},
\textsc{ANPCyT, PICT 03-08280} and \textsc{Cooperaci\'on Internacional--CONICET}. 
The second author is a research member of \textsc{CONICET} (Argentina).}
}
\end{footnotetext}

\begin{abstract}
We study isomorphisms between generalized Weyl algebras,
giving a complete answer to the quantum case of
this problem for $R=k[h]$. 

\end{abstract}

\section{Introduction}

Generalized Weyl algebras (GWA for
short), have been defined by V.~Bavula in \cite{Bavrep} and widely studied by himself and collaborators in a
series of papers (see for example \cite{Bavrep}, \cite{BJ}, \cite{BavLen}) from the point
of view of ring theory.

This kind of algebras has also been studied by T.~J.~Hodges \cite{H2} with the
name of {\em Non commutative deformations of type
$A$-Kleinian singularities}. 

Examples of GWA are, $n$-th Weyl algebras, $\mathcal{U}(\mathfrak{sl}_2)$,
primitive quotients of $\mathcal{U}(\mathfrak{sl}_2)$, its quantized versions and also the subalgebras of invariants of
these algebras under the action of finite cyclic subgroups of automorphisms.

A GWA over a field $k$ is defined as follows:
Let $R$ be an associative noetherian  $k$-algebra which is an integral domain; $A=R(\sigma, a)$, with $\sigma\in\Aut(R)$ and $a$ a non-zero central element of $R$, is the $k$-algebra
generated over $R$ by two  generators $x$ and $y$ with relations :
\begin{equation} \label{reldef2}
\begin{array}{c}
xr=\sigma(r)x, \ \forall r\in R\\
yr=\sigma^{-1}(r)y, \ \forall r\in R\\
xy=\sigma(a),\\
yx=a.
\end{array}
\end{equation}
The algebra $A$ is a noetherian domain (see \cite{Bavrep}).

Concerning the problem of isomorphism between two GWA $A=R(\sigma, a)$ and $A'=R(\sigma, a')$
raised in \cite{H2}, it is known after \cite{BJ} that:

\begin{itemize}
 
\item If $R=k[h]$ and $\sigma(h)=h-1$, then
$A\cong A' \iff \exists \eta \in k^*, \beta \in k$ such that $a(h)=\eta a'(h-\beta)$.

\item If $R=k[h, h^{-1}]$ and $\sigma(h)=qh$, $q\in k^*$ not a root of unity, then
$A\cong A' \iff \exists \eta,\beta \in k^*$, $m \in \mathbb{Z}$ and $\epsilon \in 
\{1, -1\}$ such that $a'(h)=\eta h^m a'(\beta h^{\epsilon})$.

\end{itemize}

The main result of this paper gives an answer to the quantum case of
this problem for $R=k[h]$. Namely we obtain the following Theorem.

\begin{theoreme}
Set $q\in k^*$ not a root of unity, $a_1,a_2\in k[h]$ non-constant polynomials and $\sigma(h)=qh$. 
Then there is an isomorphism $\Phi$ from
 the quantum GWA $A_q(a_1)=k[h](\sigma, a_1)$  onto $A_q(a_2)=k[h](\sigma, a_2)$
 if and only if
 there exist $\alpha,\rho\in k^*$ such that $a_2(h)=\rho a_1(\alpha h)$.
\end{theoreme}

The above theorem covers for example the cases of the quantum Weyl algebra
$A_1^q$, $k$-generated by variables $x,y$ such that $xy-qyx=1$, and the algebras
$D_{q, \lambda}$ defined in \cite{ADNagoya}, that is, the quotients of the ad-locally finite elements $F_{q}$ of $\mathcal{U}_q(\mathfrak{sl}_2)$ by the ideal generated by $(\Omega - \lambda)$.
The algebra $D_{q, \lambda}$ verifies that $\Aut(D_{q, \lambda}) \cong \Aut(F_q)$.
In this case, our result says that $D_{q, \lambda}\cong D_{q, \lambda'}$  if and only
if $\lambda^2 =\lambda'^2$.

We also solve the problem of isomorphisms of quantum Smith algebras, generalizing 
Theorem 5.4 of \cite{BJ}. This allows us to treat some cases of Witten's deformations of 
$\mathcal{U}(\mathfrak{sl}_2)$ and of conformal $\mathfrak{sl}_2$ algebras (see \cite{BR}).

Finally we give a necessary condition for Morita equivalence of GWA for $R=k[h]$. This problem
has been solved  by Hodges \cite{H1} for  $\sigma(h)=h-1$ and ${\rm deg}(a)=2$. Here we look at the quantum case for arbitrary $a$.

\medskip

The article is organized as follows:

In section \ref{sect:iso} we recall the canonical forms of GWA, for $R=k[h]$, separating
the classical case from the quantum one by means of two invariants defined by Alev and Dumas in \cite{ADfracaq}. We also rephrase Theorem 3.3 of \cite{BJ} (isomorphisms between classical 
GWA)
in terms of the roots of the defining polynomials of the algebras. We then prove a technical lemma and our main result.

In section \ref{sect:morita} we 
give a necessary condition for Morita equivalence.

We will work over an algebraically closed field $k$ of characteristic zero and all algebras will be
$k$-algebras. By ``domain'' we shall mean an integral ring, not necessarily commutative.

We consider GWA of degree one : $A=k[h](\sigma, a)$, with $\sigma\in\Aut(k[h])$ and $a\in k[h]\setminus\{0\}$, is the $k$-algebra
generated over $k[h]$ by two  generators $x$ and $y$ with relations :
\begin{equation} \label{reldef}
\begin{array}{c}
xh=\sigma(h)x,\\
yh=\sigma^{-1}(h)y,\\
xy=\sigma(a),\\
yx=a.
\end{array}
\end{equation}

Anywhere in this text ``Morita equivalence'' means ``$k$-linear Morita equivalence'', and  morphisms are morphisms of $k$-algebras.

We will suppose that $a$ is not constant, because if $a=a_0 \in k^*$, then 
$A$ is the Ore extension $k[h](x^{\pm 1}; \sigma, 0)$ whose properties have been extensively studied.

We want to thank Fran\c cois Dumas for bringing to our attention the article \cite{ADNagoya}.

\section{Isomorphism and rational equivalence of GWA.}\label{sect:iso}

\subsection{General results.}
Let us first present ``canonical" forms of  GWA of degree one and of their fraction fields.

\begin{sprop} \label{prop1}
Let $A=k[h](\sigma, a)$ be a GWA of degree one. Then necessarily $\sigma(h)=qh-h_0$ with $q\in k^*$ and $h_0\in k$, and we have
one of the three following alternatives :
\begin{enumerate}
\item  $\sigma={\rm id}$, i.e. $(q, h_0)=(1,0)$, then $A$ is commutative, and $\Frac(A)$ is the (commutative) field of rational functions $k(h,x)$ ;
\item  $q=1$ and $h_0\neq 0$, then $A$ is isomorphic to
$k[h'](\sigma_{cl},a')$, with $\sigma_{cl}(h')=h'-1$ and $a'(h')=a(h_0h')$, and $\Frac(A)$ is the first Weyl division algebra ${\mathcal D}_1(k)$ 
(``classical" case) ;
\item  $q\neq 1$, then $A$ is isomorphic to $k[h'](\sigma_q, a')$, with
$\sigma_q(h')=qh'$ and $a'(h')=a(h'-h_0/(1-q))$, and $\Frac(A)$ is the quantum skew-field $k_q(x,h'
)$ 
 (``quantum" case).
\end{enumerate}
\end{sprop}
\begin{demo}
The first case is obvious. In the second one the ``canonical" form is obtained by setting $h'=h_0^{-1}h$, and the division algebra $\Frac(A)$ is the fraction field of the
$k$-algebra generated by $x$ and $h$ because $y=a(h)x^{-1}$.
It is also the fraction field of the $k$-algebra generated  by $x$ and $y'=-h_0^{-1}hx^{-1}$, and one  easily checks that $[x,y']=1$, so this $k$-algebra  is nothing but the first Weyl
algebra (because we are in characteristic zero, any homomorphic image of $A_1(k)$ is an isomorphic image).

Suppose now that $q\neq 1$, the isomorphism $\sigma_q$
is obtained by setting $h'=h+h_0/(1-q)$, and again because $y=a(h)x^{-1}$, we see that $\Frac(A)$ is the fraction field of the $k$-algebra generated by $x$ and $h$, or by $x$ and $y'=hx-xh=((1-q)
h+h_0)x$.
One can then easily check that $xy'=qy'x$, and these two generators are algebraically independent over $k$.
\end{demo}

\medskip

The following invariants, introduced by Alev and Dumas in \cite{ADfracaq},   separate these three cases.

\begin{sdefin} \label{defeg}
For any  associative $k$-algebra $A$, one notes :
\begin{itemize}
\item $G(A)=(A^*)'\cap k^*$ the trace on $k^*$ of the derived group from the   multiplicative group $A^*$ of units in $A$;
\item $E(A)=[A,A] \cap k$ the trace on $k$ of the derived Lie algebra  from $A$.
\end{itemize}
\end{sdefin}

{\bf Remark.}
The space $E(A)$ can only be $k$  or $\{0\}$, according to the  existence or not   of two elements $x,y\in A$ such that $[x,y]=1$.

\medskip

The following proposition separates   Weyl division algebras and quantum skew-fields.

\begin{sprop} \label{egcw}
\begin{enumerate}
\item The first Weyl division algebra ${\mathcal D}_{1}(k)$ satisfies :
$$G({\mathcal D}_{1}(k))=\{1\},\ \textrm{ and }\  E({\mathcal D}_{1}(k))=k.$$
\item Set $q\in k^*$, then the quantum skew-field ${k}_{q}(u,v)$ satisfies :
$$E({k}_{q}(u,v))=\{0\},\ \textrm{ and }\
G({k}_{q}(u,v))=\ \langle q\rangle,$$
the multiplicative subgroup of $k^*$ generated by $q$.
\end{enumerate}
\end{sprop}
\begin{demo}
These are particular cases of proposition 3.9 and Th\'eor\`eme 3.10 of \cite{ADfracaq}.
\end{demo}

\begin{scorol} \label{cor1}
The Weyl algebra $A_{1}(k)$ never imbeds in the fraction field of a  \emph{noncommutative} quantum plane $k_q[u,v]$, and vice-versa. 
\end{scorol}
\begin{demo}
The proof  is direct from proposition \ref{egcw}.
\end{demo}

\medskip

{\bf Remark.}
The GWA $k[h](\sigma, a(h))$ is isomorphic to $k[h](\sigma^{-1},a(\sigma(h)))$ by $(x,y)\mapsto(y,x)$, so $k[h](\sigma_q,a(h))\cong k[h](\sigma_{q^{-1}},a(qh))$. Because two
isomorphic GWA must be rationally equivalent, and thanks to propositions \ref{prop1} and 
\ref{egcw} we may
consider the isomorphism problem for fixed $\sigma=\sigma_{cl}$ or $\sigma=\sigma_q$, $q\in k^*$.

In the classical case, the isomorphism problem is completely solved by Bavula
and Jordan in \cite{BJ} in the following way.

\begin{stheoreme} \label{isobj}
Let $A=k[h](\sigma_{cl},a_1)$ and $B=k[h](\sigma_{cl},a_2)$ be two classical
GWA, with $\sigma_{cl}(h)=h-1$, and $a_1,a_2\in k[h]$. Then $A$ and $B$ are
isomorphic if and only if $a_2(h)=\rho a_1(\epsilon h +\alpha)$, for some $\rho\in
k^*$,  $\alpha\in k$
and $\epsilon\in\{-1,1\}$.
\end{stheoreme}
\begin{demo}
This is Theorem 3.3 of \cite{BJ}.
\end{demo}

\medskip

Let us rephrase this result in terms of the roots of the polynomials $a_1$ and
$a_2$ in the  case where these are non-constant polynomials.

\begin{scorol} \label{bjroots}
Let $A=k[h](\sigma_{cl},a_1)$ and $B=k[h](\sigma_{cl},a_2)$ be two classical
GWA, with $\sigma_{cl}(h)=h-1$, and $a_1,a_2\in k[h]\setminus k$. 
Write $a_1=\prod_{i=1}^n(h-\alpha_i)$ and $a_2=\prod_{i=1}^m(h-\beta_i)$
Then the following assertions are equivalent :
\begin{enumerate}
\item $A$ and $B$ are isomorphic ;
 \item $n=m$, and there exists $\tau\in {\mathcal S}_n$  such that for all $1\leq i,j\leq n$, one has $(\a_i-\a_j)^2=(\b_{\tau(i)}-\b_{\tau(j)})^2$ ;
 \item $n=m$, and there exist $\tau\in {\mathcal S}_n$ and $\epsilon\in \{-1,1\}$ such that for all $1\leq i\leq n$, one has $\a_i-\a_1=\epsilon(\b_{\tau(i)}-\b_{\tau(1)})$. 
\end{enumerate}
\end{scorol}
\begin{demo}
It is clear that by Theorem \ref{isobj} the first point implies points 2 and 3. Clearly point 3 implies point 2.
So we only have to check that if $a_1$ and $a_2$ satisfy the second point,  then we have $a_2(h)=\rho a_1(\epsilon h +\alpha)$, with $\rho\in
k^*$,  $\alpha\in k$
and $\epsilon\in\{-1,1\}$.

Suppose first that both $a_1$ and $a_2$ only have simple roots.
By hypothesis, for each $(i,j)\in \{1,\ldots,n\}^2$ with $i\neq j$ there is $\epsilon_{i,j}\in\{-1,1\}$ such that $\a_i-\a_j=\epsilon_{i,j}(\b_{\tau(i)}-\b_{\tau(j)})$. By
Theorem \ref{isobj} $A$ is isomorphic to $k[h](\sigma_{cl},a_1(-h))$, and the roots of $a_1(-h)$ are $-\a_1,\ldots,-\a_n$.
So without loss of generality one may suppose that $\epsilon_{1,2}=1$. By reordering the roots of $a_2$, one may also suppose that $\tau=\Id$.
Let $\beta=\beta_1-\alpha_1$. We show now that the $\b_i-\beta$, $1\leq i\leq n$ are the roots of $a_1$. Since $a_1$ and $a_2$ have exactly $n$ simple roots, this will show that
$a_2(h)=\rho a_1(h+\beta)$ for some $\rho\in k^*$.

By definition
$\beta_1-\beta=\alpha_1$  is a root of $a_1$, and $\b_2-\b=\b_2-(\b_1-\a_1)=\epsilon_{1,2}(\a_2-\a_1)+\a_1=\a_2$ is a root of $a_1$ too.
Now for any $i\geq 3$, we have 
$$\a_i=\a_1-(\a_1-\a_i)=\b_1-\b-\epsilon_{1,i}(\b_1-\b_i)=\epsilon_{1,i} \b_i-\b+(1-\epsilon_{1,i})\b_1.$$
In an analogous way we show that  $\a_i=\epsilon_{2,1} \b_i-\b+(1-\epsilon_{2,i})\b_2$.
In particular 
$$(1-\epsilon_{1,i})\beta_1 + \epsilon_{1,i}\beta_i= 
(1-\epsilon_{2,i})\beta_2 + \epsilon_{2,i}\beta_i.$$


If $\epsilon_{1,i}=-1$, then 
$2\beta_1-\beta_i=(1-\epsilon_{2,i})\beta_2+
\epsilon_{2,i}\beta_i$.
This last expression equals $\beta_i$ when
$\epsilon_{2,i}=1$ and 
$2\beta_2-\beta_i$ when $\epsilon_{2,i}= -1$. In the first
case, we obtain
$\beta_1=\beta_i$, and in the second one,
$2\beta_1=2\beta_2$. Both conditions contradict our
hypothesis concerning the roots of $a_2$.

If $\epsilon_{1,i}=1$, then $\alpha_i=\beta_i- \beta$ is a
root of $a_1$, and we are done. 

In the general case,  point 2 implies that if $\alpha_i$ has multiplicity $m$ in $a_1$, then $\beta_{\tau(i)}$ has multiplicity $m$ too in $a_2$.
Now the proof follows from the case where all the roots are different.
\end{demo}



\subsection{The quantum case.}

The isomorphism problem was solved too in \cite{BJ} for GWA of type $k[h^{\pm 1}](\sigma_q,a)$ with $a\in k[h^{\pm 1}]$.

\begin{stheoreme} \label{isobjq}
Let $A=k[h^{\pm 1}](\sigma_{q},a_1)$ and $B=k[h^{\pm 1}](\sigma_{q},a_2)$ be two quantum
GWA, with $\sigma_{q}(h)=qh$, $q\in k^*$, and $a_1,a_2\in k[h^{\pm 1}]$. Then $A$ and $B$ are
isomorphic if and only if $a_2(h)=\rho h^ma_1(\alpha h^{\epsilon})$, for some $\rho,\alpha\in
k^*$,  $m\in{\mathbb Z}$,
and $\epsilon\in\{-1,1\}$.
\end{stheoreme}
\begin{demo}
This is Theorem 5.2 of \cite{BJ}.
\end{demo}

\medskip

Let us now consider the quantum GWA $k[h](\sigma_q,a)$ with $a\in k[h]$.
The preceding result cannot be used to resolve this case, because its proof relies on the fact that $h$ is invertible in $k[h^{\pm 1}]$. 
We will see in lemma \ref{Lemmad} that $h$ generates (multiplicatively) the set of normalizing elements in $k[h](\sigma_q,a)$.

\medskip

{\bf Examples} $\bullet$ For $a\in k^*$ we obtain the following localization of the quantum plane $A_q(a)=k_q[x^{\pm 1},h]$.\\
$\bullet$ For deg$(a)=1$ the algebra $A_q(a)$ is generated by $x$ and $y$, and one gets  Manin's quantum plane $k_q[x,y]$ for $a=h$,  and the quantum Weyl algebra $A_1^q(k)$
for $a=h-h_0$, with $h_0\in k^*$.\\
$\bullet$ For deg$(a)=2$ and $a(0)\neq0$, thanks to the proposition below one can suppose without loss of generality that $a(h)=-q^{-1}/(q-1)h^2+q^{-1}\lambda
(q-1)/{q'}^2 h-1/(q-1)$, for some $\lambda\in k$, with $q'\in k^*$ such that ${q'}^4=q$. Then the algebra $A_q(a)$ is a particular subalgebra of the primitive quotient $B_{q',\lambda}$ of
the quantum enveloping algebra $U_{q'}({\mathfrak{sl}}(2))$. More precisely (see \cite{ADNagoya} for details), $A_q(a)=D_{q',\lambda}$ is the canonical image in
$B_{q',\lambda}$ of ad-locally finite elements in $U_{q'}({\mathfrak{sl}}(2))$.

\medskip

Note that the (non-Laurent) polynomial version of the condition appearing in Theorem \ref{isobjq}  is obviously sufficient to obtain isomorphic GWA.

\begin{sprop} \label{prout}
Set $q\in k^*$, and let $A_q(a_1)=k[h](\sigma_{q},a_1)$ and $A_q(a_2)=k[h](\sigma_{q},a_2)$ be two quantum GWA of degree one, with $\sigma_q(h)=qh$, and $a_1, a_2\in k[h]$ such that
$a_2(h)=\rho a_1(\beta h)$, with $\rho, \beta\in k^*$.
Then $A_q(a_1)$ and $A_q(a_2)$ are isomorphic.
\end{sprop}
\begin{demo}
The isomorphism is defined by $h\mapsto \beta h$, $x\mapsto \rho^{-1}x$ and $y\mapsto y$.
\end{demo}

\medskip

We prove now that this condition is necessary. 

\medskip

{\bf Notation.} For $q\in k^*$, and $a\in k[h]$, let $A=A_q(a)$ be the GWA $k[h](\sigma_q, a)$.
We fix for all this paragraph $q\in k^*$ not a root of unity, and 
 we assume deg$(a)\geq 1$.
 
\medskip

Recall that 
 $A$ is ${\mathbb Z}$-graded by deg$(h)=0$,  deg$(x)=1$, and deg$(y)=-1$ (see \cite{Bavrep}), so that
$$A=\bigoplus_{n>0}k[h]y^n\oplus k[h]\oplus\bigoplus_{n>0}k[h]x^n.$$

The following Lemma is the generalization to arbitrary degree of  \cite{ADNagoya}, Lemma 4.5. 

\begin{slemme} \label{Lemmad}
\begin{enumerate}
\item
Suppose that $a\in k[h]$ is not a monomial. Then the normalizing elements of $A$ are the monomials $\alpha h^n$, with $\alpha \in k$ and $n\in{\mathbb N}$.
\item Suppose that $a\in k[h]$ is  a monomial. Then the normalizing elements of $A$ are the monomials $\alpha x^ph^n$, with $\alpha \in k$ and $p,n\in{\mathbb N}$, and $\alpha
y^ph^n$, with $\alpha \in k$ and $p,n\in{\mathbb N}$.
\end{enumerate}
\end{slemme}
\begin{demo}
Let $d\in A$ be a normalizing element.
Set $T_i=k[h]x^i$ if $i\geq 0$ and $T_i= k[h]y^{-i}$ if $i<0$, so that $A=\oplus_{n\in{\mathbb Z}}T_n$, and
 write $d=\sum_{i=i_1}^{i_2} \tilde d_i$, with $i_1\leq i_2\in{\mathbb Z}$, and $\tilde d_i\in T_i$.
By hypothesis, there exists $s\in A$ such that $hd=ds$. Write again $s=\sum_{j=j_1}^{j_2}\tilde s_j$, with $\tilde s_j\in T_j$.
Necessarily $hd\in \oplus_{i=i_1}^{i_2}T_i$, so we must have $s=\tilde s_0\in k[h]$.
Then  $hd=ds$ implies $h\tilde d_{i_1}=\tilde d_{i_1}\tilde s_0$ and $h\tilde d_{i_2}=\tilde d_{i_2}\tilde s_0$. Since $h\tilde d_i=q^{-i}\tilde d_ih$ and $q$ is not a root of
unity we get $i_1=i_2$, i.e. $d\in T_{i_1}$, and $\tilde s_0=q^{-i_1}h$.

Assume $i_1<0$, and put $\iota=-i_1=|i_1|$, so that $d=d(h)y^{\iota}$, and $\iota\geq 1$.
There exists $t\in A$ such that $xd=dt$. Once again the grading of $A$  implies that $t\in T_1$, and we denote  $t=t(h)x$. Then
$$xd=xd(h)y^{\iota}=d(qh)(xy)y^{\iota-1}=d(qh)a(qh)y^{\iota-1}.$$
Similarly, 
$$dt=d(h)y^{\iota}t(h)x=d(h)t(q^{-\iota}h)y^{\iota-1}a(h)=d(h)t(q^{-\iota}h)a(q^{1-\iota}h)y^{\iota-1}.$$
Now looking at the degree in $h$ of these two elements of $T_{i_1+1}$ we get  $t(h)=t_0\in k^*$ and $d(qh)a(qh)=t_0d(h)a(q^{1-\iota}h)$.
Write $a(h)=\sum_{i=n_1}^{n}a_ih^i$, and $d(h)=\sum_{i=m_1}^{m}d_ih^i$, with $0 \leq m_1\leq m$, $0 \leq n_1\leq n$ (we have $n_1<n$ if $a$ is supposed not to be a monomial), and  $a_{n_1}, a_{n}, d_{m_1}, d_{m}\in k^*$. Identifying the term of highest degree in $d(qh)a(qh)$
and $t_0d(h)a(q^{1-\iota}h)$ we get
$$q^{m}d_{m}q^na_nh^{m+n}=t_0d_{m}q^{n(1-\iota)}a_nh^{m+n},$$
 so $t_0=q^{m+n\iota}$.
Then identifying the term of lowest degree we get 
$$q^{m_1}d_{m_1}q^{n_1}a_{n_1}h^{m_1+n_1}=q^{m+n\iota}d_{m_1}q^{n_1(1-\iota)}a_{n_1}h^{m_1+n_1},$$
 so $m_1=m+(n-n_1)\iota$.
 
 Suppose first that $a$ is not a monomial. Then we get
a contradiction since $n_1<n$. This proves that $i_1\geq 0$. Symmetrizing the roles of $x$ and $y$, one gets in an analogous way $i_1\leq 0$, and finally $i_1=0$, and $d=d(h)\in k[h]$.
Then the equality $xd=tx$ becomes $d(qh)=d(h)t(h)$. So $t(h)=t_0\in k^*$, and $d(h)$ is homogeneous, i.e. it is a monomial.

Now if $a$ is a monomial we have $n_1=n$, so $m_1=m$, and $d(h)=d_mh^m$, i.e. $d=d_mh^my^{\iota}$. In the same way if $i_1> 0$ then $d=d_mh^mx^{i_1}$. At last, if $d=d(h)\in
k[h]$ is normalizing then it has to be a monomial.

To end the proof, we just have to check that the announced elements are normalizing, which is easily done since $\sigma=\sigma_q$.
\end{demo}

\bigskip

{\bf Remark.}
Suppose  $a=\alpha h^m$. Then  $A_q(a)$ is linearly
spanned by its normalizing elements, and can never be
isomorphic to a $A_q(a_2)$, with $a_2$ not a monomial in $k[h]$.

When both $a_1$ and $a_2$ are monomials, then $A_q(a_1)$
and $A_q(a_2)$ are isomorphic if and only if ${\rm deg}(a_1)= {\rm deg}(a_2)$ (this follows for instance from the computation of the Hochschild homology of $A_q(a)$, see
\cite{FSSA2}).

\medskip

This remark together with the previous results   give rise to the following Theorem.

\begin{stheoreme} \label{cnsiso}
Set $q\in k^*$ not a root of unity, and $a_1,a_2\in k[h]$ are non-constant polynomials. Then there is an isomorphism $\Phi$ from
 the quantum GWA $A_q(a_1)$  onto $A_q(a_2)$ if and only if
  there exist $\alpha,\rho\in k^*$ such that $a_2(h)=\rho a_1(\alpha h)$.
\end{stheoreme}
\begin{demo}
We saw in proposition \ref{prout} that the condition is sufficient. For the converse,
by the remark above we only have to consider the case where neither $a_1$ nor $a_2$ is a monomial. 
 Denote $h,x,y$ the generators of $A_q(a_1)$, and $h',x',y'$ the generators of $A_q(a_2)$.
The isomorphism $\Phi$ maps the normalizing elements of $A_q(a_1)$ onto the normalizing elements of $A_q(a_2)$. So one necessarily has $\Phi(h)=\alpha h'$, with $\alpha \in k^*$.
Now from the equality $xh=qhx$ in $A_q(a_1)$ we deduce $\Phi(x)\alpha h'=q\alpha h'\Phi(x)$ in $A_q(a_2)$. Just like in the proof of Lemma \ref{Lemmad}, thanks to the grading of $A_q(a_2)$
this equality implies that $\Phi(x)=d(h')x'$, with $d(h')\in k[h']$. Considering $\Phi^{-1}$ instead of $\Phi$, one gets in the same way that $\Phi^{-1}(x')=e(h)x$.
So $x=\Phi^{-1}(\Phi(x))=d(\alpha^{-1}h)e(h)x$ implies $d(h)=\beta\in k^*$, and finally $\Phi(x)=\beta x'$. Similarly one proves
$\Phi(y)=\gamma y'$, with $\gamma\in k^*$.
Then in $A_q(a_2)$ one has
$$a_1(\alpha h')=a_1(\Phi(h))=\Phi(a_1(h))=\Phi(yx)=\gamma y'\beta x'=\gamma\beta a_2(h'),$$
and we are done.
\end{demo}

\medskip

{\bf Remark.} By analogy with corollary \ref{bjroots}, one can write in the non-monomial case a multiplicative condition on the roots of $a_1$ and $a_2$ equivalent to the isomorphism of $A_q(a_1)$ and
$A_q(a_2)$.

\ex After the Theorem we can see that $D_{q,\lambda} \cong D_{q,\lambda'}$ if and only 
if $\lambda^2= \lambda'^2$.

\medskip

From the proof of this theorem one deduces easily the automorphism group of $A_q(a)$ in the following way.

\begin{stheoreme}

Suppose $a\in k[h]$ is  non-monomial. Denote $a=\sum_{i=i_0}^n a_ih^i$, with 
$a_i\in k$ and $a_{i_0}a_n\neq 0$, and define  
$p={\rm gcd}(|i-j|, \ i,j\in{\mathbb N}, i\neq j, \ a_ia_j\neq 0)$.
Then ${\rm Aut}A_q(a)\simeq ({\mathbb Z}/p{\mathbb Z})\times k^*$.

More precisely, any automorphism of $A_q(a)$ is of the form
$$x\mapsto \alpha x,\ y\mapsto \gamma^{i_0}\alpha^{-1}y;\ h\mapsto \gamma h,$$
with $\alpha\in k^*$ and $\gamma$ a $p$th-root of unity.
\end{stheoreme}
\begin{demo}
From the proof of theorem 2.2.4 we can deduce that in the non-monomial case, any automorphism $\sigma$ of $A_q(a)$ is defined by $\sigma(x)=\alpha x$,
$\sigma(y)=\beta y$, $\sigma(h)=\gamma h$, with $\alpha, \beta, \gamma \in k^*$.
Choose then $i>j$ such that $a_ia_j\neq 0$, then from $\sigma(y)\sigma(x)=a(\sigma(h))$ one gets $\gamma^{i-j}=1$. Then one deduces (for instance using Euclide algorithm) that  $\gamma^p=1$. 
Looking at the term of lowest degree in $h$ in these polynomials one gets $\beta=\gamma^{i_0}\alpha^{-1}$, where $i_0$ is the valuation of polynomial $a(h)$.
 So $\sigma(x)=\alpha x$, $\sigma(h)=\gamma h$ and
$\sigma(y)=\gamma^{i_0}\alpha^{-1}y$.

Conversely, given $\omega$ a $p$th-root of unity in $k^*$ and $\alpha\in k^*$, one checks easily that $x\mapsto \alpha x$, $y\mapsto \omega^{i_0} \alpha^{-1}$ and
$h\mapsto \omega h$ defines an isomorphism of $A_q(a)$.
\end{demo}

\bigskip

{\bf Remark.} The previous theorem  generalizes the computation of the automorphism group of $D_{q,\lambda}$ made in [2], Theorem 4.6.

\medskip

In the monomial case, the proof is rather similar to the previous one, but we have first to check that any automorphism of $A_q(k)$ must send $h$ to $\gamma h$, with $\gamma \in
k^*$.
So assume $a=a_nh^n$, with $n\geq 1$, and $a_n\in k^*$.

\begin{sprop}
Let $\sigma\in{\rm Aut}(A_q(a))$. Then there exist $\alpha, \gamma\in k^*$ such that $\sigma(x)=\alpha x$, $\sigma(h)=\gamma h$ and
$\sigma(y)=\gamma^n\alpha^{-1} y$.
\end{sprop}
\begin{demo}
Thanks to Lemma 2.2.3 we have
 $$\sigma(x)=\alpha h^{m_x}x^{n_x}y^{p_x}, \sigma(y)=\beta h^{m_y}x^{n_y}y^{p_y}, \sigma(h)=\gamma h^{m_h}x^{n_h}y^{p_h},$$
with $\alpha, \beta, \gamma \in k^*$, and $n_xp_x=n_yp_y=n_hp_h=0$.
From $\sigma(y)\sigma(x)=a_n\sigma(h)^n$ we get that $n_x=n_y=0 \Rightarrow n_h=0$, and $p_x=p_y=0\Rightarrow p_h=0$.
Because $\sigma$ is surjective we must have $n_x+n_y+n_h\neq 0$ and $p_x+p_y+p_h\neq 0$, so we only have two possibilities : $n_xp_y\neq 0$ or $n_yp_x\neq 0$.

Consider  firstly the case $n_xp_y\neq 0$, 
and assume that $p_y\geq n_x$. Then 
$$\begin{array}{c}\sigma(y)\sigma(x)=\beta\alpha h^{m_y}y^{p_y}h^{m_x}x^{n_x}=\beta\alpha q^*h^{m_x+m_y}y^{p_y-n_x}y^{n_x}x^{n_x}\\
=a_n^{n_x}\beta\alpha
q^*h^{m_x+m_y+nn_x}y^{p_y-n_x},\end{array}$$
 where $*$ represents some  integer.

So $\sigma(h)=\gamma h^{m_h}y^{p_h}$, and $a\sigma(h)^n=a\gamma^n q^* h^{nm_h}y^{np_h}$. We have now the following equalities : 
$$\left\{\begin{array}{rcl}
m_x+m_y+nn_x&=&nm_h\\
p_y-n_x=np_h \end{array}\right.$$
Now, from $yh=q^{-1}hy$ we get $\sigma(y)\sigma(h)=q^{-1}\sigma(h)\sigma(y)$, that is
$$\beta\gamma q^{-p_ym_h}h^{m_y+m_h}y^{p_y+p_h}=\gamma\beta q^{-1-p_hm_y}h^{m_h+m_y}y^{p_h+p_y}.$$
 Because $q$ is not a root of unity we get $-p_ym_h=-1-p_hm_y$.
Replacing $p_y$ by its expression above we get $-(np_h+n_x)m_h=-1-p_hm_y$, that is $p_h(nm_h-m_y)+n_xm_h=1$. Using the first equality above we finally have
$p_h(m_x+nn_x)+n_xm_h=1$. For $n\geq 2$ we must have $p_h=0$ and $n_x=m_h=1$. From $m_x+m_y+nn_x=nm_h$ we get then $m_x=m_y=0$, and from $p_y-n_x=np_h$ we get $p_y=1$.

For $n=1$ if $p_h$=0 we get the same conclusion, elsewhere we have $p_h=n_x=1$ and $m_x=m_h=0$. But then $m_x+m_y+nn_x=nm_h$ is impossible, so we get a contradiction.
Finally for any $n\geq 1$ we have $\sigma(x)=\alpha x$, $\sigma(y)=\beta y$ and $\sigma(h)=\gamma h$.

\smallskip

Assuming that $p_y<n_x$, analoguous computations lead to $p_y=n_x=1$, a contradiction.

\smallskip

The second case to consider is  $n_yp_x=0$. The reader may check that once again the same type of computations lead to expressions of the form ``sum of nonnegative
integers = -1", a contradiction.

At last, because $\sigma(y)\sigma(x)=a\sigma(h)^n$, we get $\alpha\beta=\gamma^n$.
\end{demo}

\medskip

Then one easily obtains the following.

\begin{scorol}
Suppose $a(h)=a_nh^n$ is monomial. Then ${\rm Aut}A_q(a)\simeq k^*\times k^*$.
More precisely, any automorphism of $A_q(a)$ is of the form
$$x\mapsto \alpha x,\ y\mapsto \gamma^{n}\alpha^{-1}y;\ h\mapsto \gamma h,$$
with $\alpha,\gamma \in k^*$.
\end{scorol}
\qed

The last part of this section is devoted to some application to Smith algebras.

\subsection{Smith algebras.}

Let $q\in k^*$ not a root of unity, and $\sigma=\sigma_q$.
We consider the $k$-algebra $R(f)$ freely generated by $x,y,h$ modulo relations
  \[
  xh=\sigma(h)x, \qquad
  hy=y\sigma(h), \qquad
  xy-yx=f(h).
  \]
where $f(0)=0$. This hypothesis ensures that there exists a polynomial $a \in k[h]$ such
that  $f(h)= a(\sigma(h))- a(h)$.  The algebra $R(f)$ has been studied by S.~Smith in \cite{smith}. 
Note that the polynomial $a$ is determined up to its constant term.

 Our
interest in this algebra comes from the fact that $A_q(a)$ is the quotient of $R(f)$ by the two sided ideal
generated by $\Omega= yx-a(h)$. In particular, $R(f)$ has a filtration, setting $dg(x)=dg(y)=n$ and $dg(h)=2$ (where $n$ is the degree of the polynomial $f$),
inducing a filtration on $A_q(a)$,
and it is clear that the associated graded algebra of $R(f)$ is the quantum affine space $k<x,y,h| xy=yx; xh=qhx; yh=q^{-1}hy >$.

\smallskip

{\bf Remark.} One can easily check that $\Frac(R(f))$ is the fraction field of the $k$-algebra generated over the centre $k[\Omega]$ by $x$ and $h$ satisfying $xh=qhx$.

\smallskip

We will prove the following Theorem, generalizing Theorem 5.4 of \cite{BJ}.

\begin{stheoreme}\label{smiththm}
Suppose that we have two Smith algebras, $R(f_1)$ and $R(f_2)$, with the above notations.
Then there is
 an isomorphism $\Phi$
between them if and only if there exist $\alpha\in k$, $\beta,\rho\in k^*$ such that $a_1(h)=\rho a_2(\beta h)+\alpha$.
\end{stheoreme}
\begin{demo}
It is easy to show that the centre of $R(f_i)$ is generated by $\Omega_i$.
The isomorphism $\Phi$  then sends $\Omega_1$ to $\beta \Omega_2+\alpha$, for some $\beta \in k^*$ and $\alpha\in k$.
So $\Phi$ induces an isomorphism between $A(a_1-\alpha)$ and $A(a_2)$.
We can now conclude thanks to Theorem \ref{cnsiso}.
\end{demo}

\medskip

{\bf Example.}
In his article \cite{W}, Witten introduces a 7-parameter deformation of the universal enveloping algebra ${\mathcal U}(sl_2)$. This deformation is a unital associative algebra over ${\mathbb C}$, with generators $x,y,z$ satisfying relations
$$\begin{array}{rcl}
xz&=&(\epsilon_2+\epsilon_1z)x,\\
zy&=& y(\epsilon_4+\epsilon_3z),\\
yx-\epsilon_5xy&=& \epsilon_6z^2+\epsilon_7z, \end{array}$$
where $\epsilon=(\epsilon_1,\ldots,\epsilon_7)\in{\mathbb C}^7$. This algebra is denoted ${\mathcal M}(\epsilon)$.

If
\begin{equation}\label{witten}
\epsilon_1=\epsilon_3\not\in\{0,-1,1\},\ \epsilon_2=\epsilon_4,\ \epsilon_5=1 \ \textrm{ and } \ \epsilon_7={{2\epsilon_6\epsilon_1\epsilon_2}\over{(\epsilon_1^2-1)}},
\end{equation}
then ${\mathcal M}(\epsilon)$ is isomorphic to the Smith algebra with generators $x,y$ over ${\mathbb C}[z]$, $\sigma(z)=\epsilon_2+\epsilon_1z$ and $f(z)=a(\sigma(z))-a(z)=\epsilon_6z^2+\epsilon_7z$, where $a(z)$ is defined up to a constant by 
$$a(z)={{\epsilon_6}\over{\epsilon_1^2-1}}z^2-{{\epsilon_6\epsilon_2^2}\over{\epsilon_1(\epsilon_1^2-1)}}z+a_0.$$
By using theorem \ref{smiththm}, this algebra is isomorphic to the Smith algebra with defining polynomial $a(z)=z^2-(\epsilon_2^2/\epsilon_1)z$.

If $\epsilon_1=1$, then $a(z)=(\epsilon_7/(2\epsilon_2))z^2-((\epsilon_7\epsilon_2)/2)z+a_0$, for $\epsilon_2\neq 0$, and the algebra is isomorphic to the Smith algebra with defining polynomial $\tilde a(z)=z^2-z$.
In this situation $\epsilon_6$ must be 0.

If $\epsilon_1=-1$, $\epsilon_6=0$ and $\epsilon_2\neq 0,1$ then we can choose $a(z)=z^2-z$.

Conversely, any Smith algebra with ${\rm deg}f(z)=2$ and $f(0)=0$ is isomorphic to a Witten deformation algebra ${\mathcal M}(\epsilon)$ whose parameters satisfy (\ref{witten}).

\medskip

Also, a deformation algebra has a filtration. Le Bruyn (cf. \cite{L1} and \cite{L2}) studied the algebras ${\mathcal M}(\epsilon)$ whose associated graded algebra are Auslander-regular.
He proved that they determine a 3-parameter family of Witten deformation algebras,  with defining relations:
$$\begin{array}{rcl}
xz&=&(1+\alpha z)x,\\
zy&=&y(1+\alpha z),\\
yx-\gamma xy&=&\beta z^2+z.\end{array}$$
When $\gamma=1$ this algebra is called ``conformal ${\mathfrak sl}_2$-algebra".
When $\alpha\neq\pm 1$, it is isomorphic to the Smith algebra over ${\mathbb C}[z]$ with generators $x$ and $y$,  $\sigma(z)=1+\alpha z$, and
$a(z)=(\beta/(\alpha^2-1))z^2+(\beta/(\alpha(1-\alpha^2)))z+a_0$.
We can choose, up to isomorphism, $a(z)=z^2-z/\alpha$.

\medskip

Remarks similar to those ones of the above example for down-up algebras instead of Smith algebras can be found in \cite{BR}.


\section{Morita equivalence.}\label{sect:morita}
\subsection{A necessary condition.}

The following Lemma provides a necessary condition for Morita equivalence of
two $k$-algebras.

\begin{slemme} \label{cn}
Let $A$ and $B$ be two $k$-algebras, and assume they are  noetherian domains. If they are Morita equivalent, then their fraction fields $\Frac(A)$
and $\Frac(B)$ are isomorphic.
\end{slemme}
\begin{demo}
The rings $A$ and $B$ are Morita equivalent, that is,
the category of (for example) left $A$-modules, $A-Mod$ is equivalent
to $B-Mod$.
We use now the Example 10.3.16 of \cite{Weibel}.
Let $\Sigma_A$ and $\Sigma_B$ be the collections of morphisms
$f\in A-Mod$ (resp. in $B-Mod$) that  become isomorphims 
after tensoring over $A$ (resp. $B$) by ${\rm id}_{Frac(A)}$ (resp. ${\rm id}_{Frac(B)}$).
$\Sigma_A$ and $\Sigma_B$ are multiplicative systems, and the category
$Frac(A)-Mod$ is localizing, since the canonical map
$M\to Frac(A)\otimes_A M$ belongs to $\Sigma_A$, for
all $M\in A-Mod$.

 Since $\Sigma_A\cap Frac(A)-Mod$ consists
of isomorphims, we have that $Frac(A)-Mod\cong \Sigma_A^{-1}(A-Mod)$.
Similar facts hold for $B$.
Also, $f\in \Sigma_A$ if and only if its image under the equivalence functor
is in $\Sigma_B$.

As a consequence $Frac(A)$ and $Frac(B)$ are Morita equivalent rings,
but since they are division rings,
the respective endomorphism rings are simple, and then they are isomorphic.
\end{demo}

\medskip

{\bf Remarks} $\bullet$ One may only suppose that one of the algebras is noetherian, because this is  a Morita-invariant property (but not the  property of being an integral domain).

$\bullet$ This Lemma separates up to Morita equivalence the Weyl algebra $A_n(k)$ from  algebras $S_{n,n}^{\Lambda}(k)$ of twisted differential operators on an affine quantum space
 studied in \cite{Rhhqw}, where these algebras where shown to have the same Hochschild homology and cohomology.

\medskip

Now we  apply this Lemma to  GWA.
Note that simplicity is a Morita invariant property too, and that a quantum GWA $A_q(a)=k[h](\sigma_q,a)$ is never simple because $h$ is normal and generates a non-trivial two-sided ideal.
To have a simple algebra one may localize by $h$, and consider $A'_q(a)=k[h^{\pm 1}](\sigma_q,a)$.
Then one has the following criteria for simplicity.
Note that a Laurent polynomial $a(h)$ can always be uniquely written as $a(h)=h^m\tilde a(h)$, with $m\in{\mathbb Z}$, $\tilde a(h)\in k[h]$ and $\tilde a(0)\neq 0$.
By roots of a Laurent polynomial $a$ we mean the (non-zero) roots of the associated polynomial $\tilde a$. 

\begin{sprop} \label{simple}
\begin{enumerate}
\item The classical GWA $A(a)=k[h](\sigma_{cl},a)$ is simple if and only if the polynomial $a$ only has simple roots, and for any two distinct roots $\a$ and $\beta$ of $a$, then $\a-\beta\not\in{\mathbb
Z}$.
 \item The algebra $A'_q(a)=k[h^{\pm 1}](\sigma_q,a)$ is simple if and only if $q$ is not a root of unity, the Laurent polynomial $a$ only has simple roots, and for any two distinct roots $\a$ and $\beta$ of $a$,  there exist no
 $m\in{\mathbb Z}$ such that $\a=q^m \beta$. 
\end{enumerate}
\end{sprop}
\begin{demo}
This follows from proposition 2 and corollary 2 of \cite{BaIdeals}.
\end{demo}

\medskip

It is clear that, just like in proposition \ref{prop1},  the algebra $A'_q(a)=k[h^{\pm 1}](\sigma_q,a)$ has  fraction field the quantum skew-field $k_q(u,v)$. So applying Lemma
\ref{cn} together with the results of the first section we obtain the following.

\begin{sprop} \label{cngwa}
\begin{enumerate}
\item A  classical GWA
 $A(a_1)=k[h](\sigma_{cl}, a_1)$ is never Morita equivalent to a  quantum  GWA.
\item If  $A_{q_1}(a_1)$ and $A_{q_2}(a_2)$ are two Morita equivalent quantum GWA,
then  $\langle q_2\rangle=\langle q_1\rangle$ in $k^*$.
\item If  $A'_{q_1}(a_1)$ and $A'_{q_2}(a_2)$ are two  Morita equivalent simple quantum GWA, 
then  $q_2=q_1^{\pm 1}$ in $k^*$.
\end{enumerate}
\end{sprop}
\begin{demo}
1. If $A(a_1)=k[h](\sigma_{cl}, a_1)$ is Morita equivalent to a  noetherian domain $B$, then $Frac(A(a_1))={\mathcal
D}_1(k)$ is isomorphic to $Frac(B)$ by Lemma \ref{cn},  and by corollary \ref{cor1} $Frac(B)$ cannot be a quantum skew-field, so $B$ cannot be a quantum GWA.

Points
2 and 3 are proved in an analogous way, using the invariant $G$ defined in \ref{defeg}. For the conclusion of point 3, note that in the simple case $q_1$ and $q_2$  are not  roots of
unity.
\end{demo}

\medskip

{\bf Remarks.}
$\bullet$ In the classical case, the Hochschild homology and cohomology are computed in \cite{FSSA}, giving necessary conditions for Morita equivalence linked to the degree of
the polynomials $a_1$ and $a_2$.\\
$\bullet$  As already mentioned after corollary \ref{cor1}, the GWA $k[h^{\pm 1}](\sigma, a(h))$ is isomorphic to $k[h^{\pm
1}](\sigma^{-1},a(\sigma(h)))$ by $(x,y)\mapsto(y,x)$. So in  point 3 one may conclude that up to isomorphism
we have $q_2=q_1$.

\subsection{Remarks about the classification problem.}
{\bf The classical case.}

For $a\in k[h]$, let $A(a)$ denote the GWA $k[h](\sigma_{cl},a)$.

$\bullet$ If deg$(a)=1$, then $A(a)$ is isomorphic to the Weyl algebra $A_1(k)$.

$\bullet$ If deg$(a)=2$, then $A(a)$ is a prime quotient of the enveloping algebra ${\mathcal U}({ sl}_2)$, and the classification up to Morita equivalence is completely
solved by Hodges in \cite{H1} by means of the following Theorem.
\begin{stheoreme} \label{hodges}
Set $a_1,a_2\in k[h]$ two polynomials of degree 2 with distinct roots respectively $\a_1\neq\a_2$ and $\a'_1\neq\a'_2$. Then the GWA $A(a_1)$ and $A(a_2)$ are Morita-equivalent if and only if 
there exist $\epsilon\in\{-1,1\}$ and $m\in{\mathbb Z}$ such that $\a_1-\a_2=\epsilon(\a'_1-\a'_2)+m$.
\end{stheoreme}
\begin{demo}
This is just a rephrasing in terms of roots of the polynomials of \cite{H1}, Theorem 5.
\end{demo}

\medskip

$\bullet$ If deg$(a)\geq 3$, we can prove the following sufficient condition.
\begin{sprop} \label{prophodges}
Set $a_1,a_2\in k[h]$ two polynomials of degree $n\geq 3$ satisfying the simplicity criterium of proposition \ref{simple}(1), with distinct roots respectively $\a_1,\ldots,\a_n$ and $\a'_1,\ldots,\a'_n$. Suppose that
there exist $\tau\in{\mathcal S}_n$ and $(m_1,\ldots,m_n)\in{\mathbb Z}^n$ such that $\a_i=\a'_{\tau(i)}+m_i$ for all $1\leq i\leq n$. Then the GWA $A(a_1)$ and
$A(a_2)$ are Morita-equivalent.
\end{sprop}
\begin{demo}
This is a direct consequence of Lemma 2.4 and Theorem 2.3 of \cite{H2}.
\end{demo}

\medskip

This condition is exactly the analogous in degree $n$ of the condition of Theorem \ref{hodges}. So, as already noted by Hodges,  it would be interesting to know a necessary condition in any degree.
Since the condition appearing in the above proposition is written in terms of the roots of the polynomials, so one may compare it  with corollary \ref{bjroots}.

\medskip

{\bf The quantum case.}

\medskip

We consider the non root of unity case.
At this point, we have some necessary conditions, first by proposition \ref{cngwa} for $A_{q_1}(a_1)$ and $A_{q_2}(a_2)$ to be Morita equivalent we must have $q_2=q_1^{\pm 1}$.
Moreover the degree of the polynomial used to define a GWA is
linked to its  Hochschild homology, computed by Farinati, Solotar and
Su\'arez-\'Alvarez in \cite{FSSA2}. 
Then if $A_q(a_1)$ and $A_q(a_2)$  are Morita equivalent one must also have that
${\rm deg}(a_1)={\rm deg}(a_2)$ and the same for the  polynomials ${\rm gcd}(a_1,a'_1)$ and ${\rm gcd}(a_2,a'_2)$.
We can say too that the sufficient conditions mentioned in proposition \ref{prophodges}  are still true in the quantum case, replacing in Hodges'
proofs (\cite{H2}, Lemma 2.4 and Theorem 2.3) $\sigma_{cl}$ by $\sigma_q$.
Then the condition $\a_i=\a'_{\tau(i)}+m_i$ becomes $\a_i=q^{m_i}\a'_{\tau(i)}$.

\medskip


\end{document}